\documentclass[11pt]{article}
\usepackage{amssymb}

\newenvironment{proof}{\begin{trivlist}\item[]{\it
Proof.}}{\hfill$\square$\end{trivlist}}

\newenvironment{proofofthm}[1]{\noindent{\it Proof of Theorem
#1}}{\hfill$\square$\\\mbox{}}

\newtheorem{theorem}{Theorem}[section]

\newtheorem{definition}[theorem]{Definition}
\newtheorem{lemma}[theorem]{Lemma}
\newtheorem{proposition}[theorem]{Proposition}

\def\flip{\phi}
\def\mc{{\mathbb{C}}}
\def\id{{\mathrm{id}}}
\def\qdet{{\mathrm{det}}_q}
\def\cala{{\cal{A}}}
\def\lqm{{\cal{L}}_q(M)}
\def\lkqm{{\cal{L}}_{k,q}(M)}
\def\calf{{\cal{F}}}
\def\qtr{{\mathrm{Tr}}_q}

\def\ad{{\mathrm{Ad}}}
\def\orbxi{O_{\xi}}
\def\orbclxi{O_{\xi}^{{\mathrm{cl}}}}
\def\funkqorb{\calf_{\mc(q),q}(\orbclxi)}
\def\mk{{\mathbb{K}}}
\def\funqg{\calf_q(G)}
\def\funkg{\calf_{k,q}(G)}
\def\funkkg{\calf_{\mk,q}(G)}
\def\funkqg{\calf_{\mc(q),q}(G)}
\def\funkm{\calf_{k,q}(M)}
\def\funkkm{\calf_{\mk,q}(M)}
\def\funkqm{\calf_{\mc(q),q}(M)}
\def\fung{\calf(G)}
\def\funqorb{\calf_q(\orbclxi)}
\def\funqm{\calf_q(M)}
\def\funqn{\calf_q(N)}
\def\evxi{{\mathrm{ev}}_{\xi}}
\def\evhatxi{\widehat{{\mathrm{ev}}}_{\xi}}
\def\jordxi{J(\xi)}
\def\idealxi{{\cal{I}}_{\xi}}
\def\ridealxi{{\cal{J}}_{\xi}}
\def\nilp{N}
\def\qsphere{{\cal{O}}(X_{q,\alpha,\beta})}
\def\image{{\mathrm{im}}}
\def\lcqm{{\cal{L}}_{\mc,q}(M)}
\def\lcqtd{{\cal{L}}_{\mc,q}^{t,d}}

\begin{document}

\title{A quantum homogeneous space of nilpotent matrices}
\author{M. Domokos
\thanks{Partially supported by OTKA No. T034530, T046378, and the Bolyai Fellowship.}
\\ 
\\ 
{\small R\'enyi Institute of Mathematics, Hungarian Academy of 
Sciences,} 
\\ {\small P.O. Box 127, 1364 Budapest, Hungary,} 
{\small E-mail: domokos@renyi.hu } }

\date{}
\maketitle

\begin{abstract} 
A quantum deformation of the adjoint action of the special linear group 
on the variety of nilpotent matrices is introduced. 
New non-embedded quantum homogeneous spaces are obtained 
related to certain maximal coadjoint orbits, 
and known quantum homogeneous spaces are revisited. 
\end{abstract}
\medskip
\noindent MSC: 16W35; 20G42; 
17B37; 81R50 

\noindent Keywords: coadjoint orbit;  quantum group; 
quantum homogeneous space; deformation; nilpotent matrices; 
comodule algebra

\bigskip 


\section{Introduction}\label{sec:intro}

The general linear group $GL(n,\mc)$ 
acts (from the right) 
on the space $M=M(n,\mc)$ of $n\times n$ matrices with complex entries 
by conjugation: 
\[\ad:M\times GL(n,\mc)\to M,\qquad (\xi,g)\mapsto g^{-1}\xi g.\] 
This is an instance of the adjoint (coadjoint) action of a reductive Lie group 
on (the dual of) its Lie algebra. For the relevance of coadjoint orbits 
in representation theory and mathematical physics see 
\cite{kirillov}. 

Throughout this paper $G$ stands for the special linear group 
$SL(n,\mc)$. Note that two points of $M$ belong to the same 
$GL(n,\mc)$--orbit if and only if they belong to the same $G$--orbit.  
For the sake of technical convenience we shall restrict $\ad$ to $G$,  
and look for quantum deformations of the orbits. 
The precise formulation of the problem can be given in the language of affine algebraic geometry. We shall write $\calf(V)$ for the coordinate ring of an affine algebraic variety $V$. The action $\ad$ is encoded in its comorphism  
\[\ad^*:\calf(M)\to\calf(M)\otimes_{\mc}\fung,\]
which (in addition to being an algebra homomorphism) is a right 
coaction of the Hopf algebra $\fung$ on the algebra $\calf(M)$. 
Take $\xi\in M$ and its $G$--orbit $\orbxi$ in $M$. 
Recall that 
$\orbxi$ is closed if and only if $\xi$ is semisimple (i.e. diagonalizable). 
The Zariski closure $\orbclxi$ is an affine $G$--subvariety of $M$; 
equivalently, its 
coordinate ring $\calf(\orbclxi)$ is a quotient $\fung$--comodule algebra of $\calf(M)$. In other words, we have a commutative diagram 
 \[\begin{array}{ccc} 
\calf(M) &\stackrel{\ad^*}\longrightarrow &\calf(M)\otimes_{\mc}\fung\\
\downarrow & &\downarrow \\
\calf(\orbclxi)&\longrightarrow &\calf(\orbclxi)\otimes_{\mc}\fung
\end{array}
\] 
of algebra homomorphisms, where the vertical arrows are induced by restricting functions on $M$ to $\orbclxi$. 

Our base ring to define quantum algebras will be 
$\mk=\mc[q]_{(q-1)}$, the localization of the polynomial ring $\mc[q]$ in the indeterminate $q$ at the maximal ideal $(q-1)$. In other words, $\mk$ is the subring of the field 
$\mc(q)$ of rational functions, consisting of the functions with no pole at $1$. Denote by $\funqg$ the quantum coordinate ring of $G$; 
this is a $\mk$--Hopf algebra defined in terms of generators and relations 
(see Section~\ref{sec:main}).  
Moreover, $\funqg$ is a deformation of 
$\fung$ in the following sense: 
$\funqg$ is a free $\mk$--module, and its quotient 
Hopf algebra modulo the ideal 
$(q-1)\funqg$ is isomorphic to $\fung$. 
(It is more standard to work over the ring $\mc[[h]]$ of formal power series when studying deformations of algebras, see for example \cite{bfgp}. 
However, all our constructions can be performed over $\mk$, 
and this choice of the base ring allows to relate our objects naturally to 
their classical counterpart in two different ways. 
The specialization $q\mapsto 1$ makes sense on the one hand, and on the other hand, it is also possible to extend scalars to $\mc(q)$ 
and formulate the results in a representation theoretic flavour.) 

\begin{definition}\label{def:qdef}
{\rm By a {\it quantum deformation of an affine $G$--variety} $V$ we mean a (right) 
$\funqg$--comodule algebra $\calf_q(V)$ with the following properties:
\begin{itemize}
\item[(i)] $\calf_q(V)$ is a free $\mk$--module; 
\item[(ii)] $\calf_q(V)/(q-1)\calf_q(V)\cong \calf(V)$;
\item[(iii)] The diagram 
\[\begin{array}{ccc} 
\calf_q(V) &\stackrel{\varphi}\longrightarrow &\calf_q(V)\otimes_{\mk}\funqg\\
\downarrow & &\downarrow \\
\calf(V)&\longrightarrow &\calf(V)\otimes_{\mc}\fung
\end{array}
\] 
commutes, where the horizontal arrows are the coaction maps, and the vertical arrows are the natural surjections with kernel generated by $(q-1)$. 
When $V$ contains a dense $G$--orbit, we say that 
$(\calf_q(V),\varphi)$ is a 
{\it quantum homogeneous $G$--space}. 
\end{itemize}
}
\end{definition}

Denote by $P$ the subset of $M$ consisting of the matrices $\xi$ 
which have both of the following two properties: 
\begin{itemize}
\item[(i)] all non-zero eigenvalues of $\xi$ have multiplicity one; 
\item[(ii)] $\xi$ has only one nilpotent Jordan block. 
\end{itemize} 
Note that for $\xi\in P$ the $G$--orbit of $\xi$ is maximal; that is, it has 
maximal possible dimension. 
For all $\xi\in P$ we shall construct an $\funqg$--comodule algebra $\funqorb$, which is a quantum deformation of $\orbclxi$ in the above sense. 
In particular, we obtain quantum homogeneous spaces of nilpotent matrices. 
The algebra $\funqorb$ is defined as a quotient of the reflection equation algebra $\lqm$ modulo an ideal generated by adjoint invariants. 

Quantizations of coadjoint orbits of $G$ given as quotients of the reflection equation algebra were considered in several papers, see for example 
\cite{dm1}, \cite{dm2}, \cite{dm3}, \cite{gur-sap}. 
The case of semisimple orbits is settled in \cite{dm3}. 
Non-semisimple orbits of some special type are dealt with in 
\cite{dm1}, \cite{dm2}; all of them are orbits of matrices with at most two eigenvalues, and with Jordan blocks of size one or two. 
The main reason that the methods of \cite{dm1}, \cite{dm2} are restricted to 
such orbits is that the reflection equation algebra has few characters. 
In other words, all complex solutions of the reflection equations are matrices 
of this special type by \cite{mudrov}.  

There is a passage between the reflection equation algebra  
and the quantized coordinate ring $\funqm$ developed in \cite{bm}, 
\cite{majid}. 
Exploiting this method, 
formulae from \cite{dl} provide a candidate for the ideal in the reflection equation algebra that should define $\funqorb$ for $\xi\in P$. The main  technical problem is to prove that the resulting quotient 
is a free $\mk$--module. This is handled by applying 
results of \cite{dfl}, that are based on the fact that 
sufficiently many characters of $\funqm$ are available to hit the orbit of each $\xi\in P$. 


\section{Main results}\label{sec:main}

First we recall two types of non-commutative algebras 
associated with the R-matrix of the vector representation of the 
Drinfeld-Jimbo algebra $U_q({\mathfrak{sl}}_n)$ (cf. \cite{drinfeld}, 
\cite{jimbo}). 
Let $k$ be a commutative integral domain, $q$ an invertible element in $k$, 
and $R\in{\mathrm{End}}_k(k^n\otimes k^n)$ 
given by 
\[
R^{is}_{jt}=\left\{
\begin{array}{cc}
q, &\mbox{ if }i=j=s=t;\\
1, &\mbox{ if }i=j,s=t,i\neq s;\\
q-q^{-1}, &\mbox{ if }i>j,i=t,j=s; \\
0, &\mbox{ otherwise} 
\end{array}\right.
\]
(see for example Section 9.2 of \cite{ks}). 
Set $\widehat{R}:=\flip\circ R$, where $\flip$ is the flip endomorphism 
$a\otimes b\mapsto b\otimes a$. 

The {\it quantum coordinate ring of 
$n\times n$ matrices} is the 
$k$--algebra $\funkm$ generated by the entries of $X=[x^i_j]_{i,j=1}^n$, 
subject to the relations 
\[\widehat{R} (X\otimes I)(I\otimes X) = 
(X\otimes I)(I\otimes X)\widehat{R}\] 
(with the notation of \cite{rtf}). 
It has a bialgebra structure: 
the comultiplication and the counit maps are given on the generators 
by $\Delta(x^i_j)=x^i_s\otimes x^s_j$ (we use the convention of summing 
over repeated indices), and 
$\varepsilon(x^i_j)=\delta^i_j$. 
The {\it quantum determinant} 
$\qdet=\sum_{\sigma\in{\mathrm{Sym}}_n}
(-q)^{{\mathrm{length}}(\sigma)}
x^1_{\sigma(1)}x^2_{\sigma(2)}\cdots x^n_{\sigma(n)}$ 
is a central group-like element 
in $\funkm$. The quotient bialgebra of $\funkm$ modulo the ideal generated by 
$(\qdet-1)$ has an antipode $S$ that makes it a $k$--Hopf algebra $\funkg$, 
called the 
{\it quantum coordinate ring of the special linear group}. 
We shall denote by $t^i_j$ the image of $x^i_j$ under 
the natural homomorphism $\pi:\funkm\to\funkg$.  
The right adjoint coaction of $\funkg$ on itself can be lifted to a coaction 
$\beta:\funkm\to\funkm\otimes\funkg$:  
for $x\in\funkm$ with $\Delta^{(2)}(x)=\sum x_1\otimes x_2\otimes x_3$, 
we have $\beta(x)=\sum x_2\otimes S(\pi(x_1))\pi(x_3)$. 
We shall be interested in these algebras in the case when $(k,q)$ 
is $(\mk,q)$ or $(\mc(q),q)$. 
Moreover, to simplify notation we write 
$\funqm$ and $\funqg$ instead of $\funkkm$ and $\funkkg$. 

The second algebra $\lkqm$ is  
the $k$--algebra generated by the entries of 
$L=[l^i_j]_{i,j=1}^n$, subject to the relations 
\begin{equation}\label{eq:reflection}
\widehat{R} (I\otimes L)\widehat{R} (I\otimes L) = 
(I\otimes L)\widehat{R}(I\otimes L)\widehat{R}.
\end{equation} 
This algebra is called the {\it reflection equation algebra} in 
\cite{kulish-sk}, and the algebra of {\it braided matrices} in \cite{majid2}. 
An important property of $\lkqm$ is that the map 
$l^i_j\mapsto l^a_b\otimes S(t^i_a)t^b_j$, $i,j=1,\ldots,n$, 
extends to an algebra homomorphism 
$\beta:\lkqm\to\lkqm\otimes_k\funkg$, 
yielding a right coaction of $\funkg$ on $\lkqm$. 
In fact, the $\funqg$--comodule algebra $\lqm$ is a quantum deformation 
(in the sense of Definition~\ref{def:qdef}) 
of the $SL(n,\mc)$--variety $M$ endowed with the adjoint action, 
where $\lqm$ stands for ${\cal{L}}_{\mk,q}(M)$.

Both $\funqm$ and $\lqm$ are graded, the generators having degree one. 
A crucial fact for us is that there is a $\mk$--module isomorphism 
$\Phi:\funqm\to\lqm$, intertwining the right $\funqg$--coactions 
$\beta$ on $\funqm$ and $\lqm$. 
This is explained in Section 7.4 of \cite{majid}, explicit formulae are given in (7.37); see also Section 10.3 in \cite{ks}.  
In particular, $\Phi(x^i_j)=l^i_j$ for all $i,j$, and $\Phi$ is homogeneous. 
We shall denote by $\Psi:\lqm\to\funqm$ the inverse $\mk$--module 
map to $\Phi$. 

The {\it quantum trace} of an arbitrary $n\times n$ matrix $M$ with 
entries $m^i_j$ in a $\mk$--module is defined as 
\[\qtr(M)=\sum_{i=1}^n q^{n+1-2i}m^i_i.\] 
The elements $\qtr(L^i)$, $i=1,2,\ldots$, are coinvariants for the 
coaction of $\funqg$ on $\lqm$; that is, 
$\beta(\qtr(L^i))=\qtr(L^i)\otimes 1$. 
Moreover, they are central in $\lqm$, see for example 
Theorem 10.3.8 in \cite{majid}. 

We need to recall the basic $\funqg$--coinvariants 
$\tau_1,\ldots,\tau_n$ in $\funqm$ 
with respect to the coaction $\beta$, introduced in \cite{dl}. 
The first of them, $\tau_1$, is the quantum trace $\qtr(X)$, 
the last of them, $\tau_n$ is the quantum determinant 
$\qdet$. For a $d$--element subset $I=\{i_1,\ldots,i_d\}$ of 
$\{1,\ldots,n\}$, the submatrix of $X$ consisting of the entries with 
row and column indices belonging to $I$ is a $d\times d$ quantum matrix, 
its quantum determinant will be denoted by $\qdet(I,I)$. 
For $d=1,\ldots,n$, we have 
\[\tau_d=\sum_{I=\{i_1,\ldots,i_d\}}
q^{(d(n+1)-2\sum_{i\in I}i)} \qdet(I,I),\]  
where the summation ranges over the $d$--element subsets of 
$\{1,\ldots,n\}$.  

Now take an $n\times n$ complex matrix $\xi$ from the set 
$P$ defined in Section~\ref{sec:intro}. 
In the $G$--orbit of $\xi$ choose  
a block diagonal Jordan normal form 
\[\jordxi={\mathrm{diag}}(J_r,\lambda_1,\ldots,\lambda_{n-r})
=\left(\begin{array}{cccccc}0 & 1 &  &  &  &  \\ & 0 & 1 &  &  &  \\ &  & 0 &  &  &  \\ &  &  & \lambda_1 &  &  \\ &  &  &  & \ddots&  \\ &  &  &  &  & \lambda_{n-r}\end{array}\right),\]
where $J_r$ 
is an $r\times r$ nilpotent Jordan block ($r\in\{0,1,\ldots,n\}$),  
and $\lambda_1,\ldots,\lambda_{n-r}$ are pairwise different non-zero 
eigenvalues of $\xi$. The map $X\mapsto \jordxi$ extends to a  $\mk$--algebra 
homomorphism $\evxi:\funqm\to\mk$; we write $\tau_d(\xi)$ for the image of $\tau_d$ under $\evxi$. 
Consider the elements 
$\Phi(\tau_d-\tau_d(\xi))=\Phi(\tau_d)-\tau_d(\xi)\in\lqm$, for 
$d=1,\ldots,n$. 
They are $\funqg$--coinvariants with respect to the coaction $\beta$, 
hence they are central  in $\lqm$ by the results in \cite{majid}. 
Write $\idealxi$ for the ideal of $\lqm$ generated by them. 
Since this ideal is generated by coinvariants, it is a subcomodule 
with respect to the coaction $\beta$. 
We define $\funqorb$ as the quotient $\funqg$--comodule algebra 
\[\funqorb=\lqm/\idealxi\] 
of $\lqm$. 
We shall keep the notation $\beta$ for the coaction of $\funqg$ on $\funqorb$. 

\begin{theorem}\label{thm:main1} 
Let $\xi$ be an $n\times n$ complex matrix from $P$. 
Then the $\funqg$--comodule algebra $\funqorb$ is a quantum deformation 
(in the sense of Definition~\ref{def:qdef}) 
of $\orbclxi$, the closure  of the $G$--orbit of $\xi$ in $M$. 
\end{theorem} 

Thus we obtained some new quantum homogeneous spaces. 
A particularly interesting extreme case is when 
$\xi=J_n$, the nilpotent $n\times n$ Jordan block.  
Then $\orbclxi$ is the set of all nilpotent $n\times n$ matrices,  
write $\nilp$ for this affine variety. It is the nullcone for the coadjoint 
action of $SL(n,\mc)$ on the dual of its Lie algebra, a relevant object both from the point of view of invariant theory and representation theory, see for example 
\cite{kostant}, \cite{kraft}. 

\begin{theorem}\label{thm:main2} 
Let $\funqn$ denote the quotient $\funqg$--comodule algebra of $\lqm$ modulo 
the ideal generated by $\qtr(L^d)$, $d=1,\ldots,n$. 
Then $\funqn$ is a quantum deformation of the affine $G$--variety $\nilp$ of 
nilpotent $n\times n$ complex matrices. 
\end{theorem}


\section{Proofs} \label{sec:proofs}

\begin{lemma}\label{lemma:qslfree}
Every $\mk$--submodule of $\funqg$ is free. 
\end{lemma}
\begin{proof} 
It is well known that every submodule of a free module over a principal ideal 
domain is free, see for example the remark after Theorem 5.3 in Chapter I of 
\cite{em}. Therefore it is sufficient to prove that $\funqg$ is a free 
$\mk$--module. 

We write $\funkg^d$ for the $k$-submodule of $\funkg$ spanned by the products of 
degree at most $d$ in the generators $t^i_j$. 
Then $\funkg^d$ 
is the sum of the spaces of matrix elements of the tensor powers of exponent $\leq d$ 
of the fundamental corepresentation $e_j\mapsto e_i\otimes t^i_j$ 
$(j=1,\ldots,n)$ of $\funkg$ (cf. Section 11.2.3 in \cite{ks}). 
Since $\funkqg$ and $\fung$ have essentially the same corepresentation 
theory (see for example \cite{hayashi} or \cite{nym}), we have the equality 
$\dim_{\mc(q)}(\funkqg^d)=\dim_{\mc}(\fung^d)$ for $d=0,1,2,\ldots$. 

By definition $\mk$ is a subring of $\mc(q)$; write $\cala$ for the 
$\mk$--subalgebra of $\funkqg$ generated by the $t^i_j$. 
Since these generators satisfy the defining relations of $\funqg$, the algebra 
$\cala$ is a homomorphic image of $\funqg$ (actually, they are isomorphic 
as we shall see below). 
Now we make use of the basis of the quantum coordinate ring of $n\times n$ 
matrices developed in \cite{gl}. It consists of certain 
(called {\it preferred} in loc. cit.) 
products of quantum minors. Although \cite{gl} works over a field, 
the arguments obviously show that as a $\mk$--module, 
$\funqg$ is spanned by the preferred products of quantum minors of 
size $\leq n-1$, and that the images of these elements in $\cala$ span 
$\funkqg$ as a $\mc(q)$--vector space. From the classical theory we know 
that the number of preferred products of (quantum) minors of 
size $\leq n-1$ and of total degree $\leq d$ is $\dim_{\mc}(\fung^d)$. 
Hence by the above dimension equality we obtain both that the homomorphism 
$\funqg\to\cala$ is an isomorphism and that $\funqg$ is a free $\mk$--module 
with basis consisting of the preferred products of quantum minors of size 
$\leq n-1$. 
\end{proof}

In particular, we may identify $\funqg$ with the $\mk$--subalgebra 
of $\funkqg$ generated by the elements $t^i_j$ $(i,j=1,\ldots,n)$, and 
we shall make this identification for the rest of the paper. 
Similarly, $\funqm$ is identified with a $\mk$--subalgebra of 
$\funkqm$ in the obvious way. 
Throughout this Section we shall use also the following convention. 
Given a $\mk$--submodule $W$ of $\funqm$, the symbol $\widehat{W}$ stands for 
the $\mc(q)$--vector subspace of $\funkqm$ spanned by $W$. 
Moreover, for any $\mk$--module $W$ 
we write $\overline{W}$ for the $\mc$--vector space 
obtained by taking the quotient of $W$ modulo the submodule 
$(q-1)W$, and write $\overline{w}\in\overline{W}$ for the image 
of $w\in W$ under the natural surjection. 
Note that by the defining relations of $\funqg$ it is clear that 
$\funqg/(q-1)\funqg$ is 
isomorphic to $\fung$, with an isomorphism mapping the generators 
$t^i_j$ to the corresponding coordinate functions on $G$. 

\begin{lemma}\label{lemma:right-ideal} 
The right ideal $\ridealxi$ 
of $\funqm$ generated by 
$\tau_1-\tau_1(\xi),\ldots,\tau_n-\tau_n(\xi)$
is mapped by $\Phi$ onto 
the ideal $\idealxi$ in $\lqm$. 
\end{lemma} 

\begin{proof} 
Note that $\Phi$ is not an algebra homomorphism. 
However, if $\beta(a)=a\otimes 1$ for some $a\in\funqm$, 
then $\Phi(ab)=\Phi(a)\Phi(b)$ in $\lqm$ for all $b\in\funqm$. 
This follows from the formula of Theorem 7.4.1 in \cite{majid}
(see also the formula of Proposition 34 (i) in Section 10.3.2 of 
\cite{ks}), 
which expresses the multiplication in $\lqm$ in terms of the multiplication 
in $\funqm$. Therefore, $\Phi(\ridealxi)$ is the right ideal in $\lqm$ 
generated by $\Phi(\tau_1-\tau_1(\xi)),\ldots,\Phi(\tau_n-\tau_n(\xi))$. 
This latter right ideal is a two-sided ideal, since it is generated by 
central elements. 
\end{proof}

\begin{lemma}\label{lemma:1} 
We have the equality \quad$\widehat{\ridealxi}\cap \funqm=\ridealxi$ 
\quad for all $\xi\in P$. 
\end{lemma} 

\begin{proof} 
Take sets of monomials $\Gamma_i\subset \funqm$ $(i=1,\ldots,n)$ 
in the variables $x^s_t$ such that 
\[\overline{\Lambda}=\bigcup_{i=1}^n\{\overline{(\tau_i-\tau_i(\xi))w}
\mid w\in\Gamma_i\}\] 
is a $\mc$--basis of $\overline{\ridealxi}$ 
compatible with the filtration by degree. 
Using the corepresentation theory of $\funkqg$, 
classical results of Kostant \cite{kostant}, and a commutative algebra lemma 
(cf. Lemma 5.2 in \cite{dfl}), it was proved in \cite{dfl} 
(see the last line of the proof of Proposition 5.3 in loc. cit.)
that 
\[\Lambda=\bigcup_{i=1}^n\{(\tau_i-\tau_i(\xi))w
\mid w\in\Gamma_i\}\] 
is a $\mc(q)$--basis of $\widehat{\ridealxi}$. 
An arbitrary element $f$ in the intersection $\widehat{\ridealxi}\cap \funqm$ 
can be uniquely written as a finite sum 
$f=\sum_{\lambda\in\Lambda}a_{\lambda}\lambda$, 
where the coefficients $a_{\lambda}$ are taken from $\mc(q)$. 
Assume that for some $\mu\in\Lambda$, the coefficient $a_{\mu}$ is not 
contained in $\mk$. Then there is a positive integer $d$ such that 
$(q-1)^da_{\lambda}\in\mk$ for all $\lambda\in\Lambda$, and 
$(q-1)^da_{\mu}$ is not contained in $(q-1)\mk$. 
Reducing the equality 
$(q-1)^df=\sum_{\lambda\in\Lambda}(q-1)^da_{\lambda}\lambda$ 
in $\funqm$ modulo the ideal $(q-1)$ we obtain that a non-trivial 
$\mc$--linear combination of the elements of $\overline{\Lambda}$ is zero. 
This contradiction shows that all the $a_{\lambda}$ are contained in $\mk$, 
hence $f$ is contained in $\ridealxi$. 
\end{proof}

\begin{proofofthm}~\ref{thm:main1}.
The $SL(n,\mc)$--orbit of $\xi\in P$ is maximal (i.e. is not contained in the closure of another orbit). 
Therefore the vanishing ideal of $\orbclxi\subset M$ in 
$\calf(M)=\lqm/(q-1)\lqm$ 
is $\overline{\idealxi}$ 
by classical results of \cite{kostant}. 
Thus $\funqorb$ satisfies (ii) and (iii) from Definition~\ref{def:qdef} by construction. 
The only thing we have to prove is that (i) holds as well, namely that 
$\funqorb$ is a free $\mk$--module. By Lemma~\ref{lemma:right-ideal}, 
the $\mk$--module isomorphism $\Psi:\lqm\to\funqm$ 
induces a $\mk$--module isomorphism between $\funqorb$ and 
$\funqm/\ridealxi$. 
Denote by $\evhatxi:\funkqm\to\mc(q)$ the $\mc(q)$--linear extension 
of $\evxi:\funqm\to\mk$ defined in Section~\ref{sec:main}, 
and denote by $\widehat{\beta}$ the $\mc(q)$--linear extension of 
the coaction 
$\beta:\funqm\to\funqm\otimes_{\mk}\funqg$ to $\funkqm$. 
Consider the $\mc(q)$--linear map 
$\widehat{\gamma}=(\evhatxi\otimes\id)\circ\widehat{\beta}
:\funkqm\to\funkqg$, and its restriction $\gamma$ to $\funqm$. 
The kernel of $\widehat{\gamma}$ is $\widehat{\ridealxi}$ by 
Theorem 5.4 of \cite{dfl}. (To be more precise, the statement of 
\cite{dfl} is about $\calf_{\mc(q),q}(GL(n,\mc))$, instead of $\funkqg$. 
This does not make any difference, as the proof of 
the cited result works word by word if we replace $\calf_{\mc(q),q}(GL(n,\mc))$ 
by $\funkqg$; alternatively, one can refer to the results of the last Section of 
\cite{dl2}.) It follows by Lemma~\ref{lemma:1} 
that $\ker(\gamma)=\ridealxi$. The image of 
the $\mk$--module homomorphism $\gamma$ is contained in 
$\funqg$. Therefore $\funqm/\ridealxi$ is isomorphic to a 
$\mk$--submodule of $\funqg$, hence is free by 
Lemma~\ref{lemma:qslfree}. 
\end{proofofthm}

\begin{proofofthm}~\ref{thm:main2} 
The elements $\Phi(\tau_1),\ldots,\Phi(\tau_n)$ 
generate the same ideal in $\lqm$ as the 
elements $\qtr(L),\qtr(L^2),\ldots,\qtr(L^n)$ 
by Lemma~\ref{lemma:newton} below.  
\end{proofofthm}

\begin{lemma}\label{lemma:newton} 
The elements $\Phi(\tau_1),\ldots,\Phi(\tau_n)$ generate the same 
$\mk$-subalgebra in $\lqm$ as the elements 
$\qtr(L),\ldots,\qtr(L^n)$. 
\end{lemma}

\begin{proof} 
Denote by $A$, $B$, and $C$, respectively, 
the $\mk$--subalgebra of $\lqm$ generated 
by the first set, the second set, and their union, respectively. 
These are graded subalgebras of $\lqm$. Their degree $d$ homogeneous 
components $A_d$, $B_d$, $C_d$ are finitely generated $\mk$--modules. 
By the classical Newton formulae relating the characteristic coefficients 
of a matrix with the traces of its powers, 
their images $\overline{A}_d$, $\overline{B}_d$, 
$\overline{C}_d$ 
in $\calf(M)$ coincide. 
Applying Nakayama's Lemma for finitely generated modules over the local 
principal ideal domain 
$\mk$ we obtain $C_d=A_d$ and $C_d=B_d$. 
This holds for all $d$, hence $A=C=B$. 
\end{proof} 

Note that Lemma~\ref{lemma:newton} implies the existence of Newton formulae 
relating $\qtr(L),\ldots,\qtr(L^n)$ and $\Phi(\tau_1),\ldots,\Phi(\tau_n)$. 
A version of Newton formulae in $\lqm$ is given in \cite{iop}. 


\section{Non-embedded quantum homogeneous spaces}\label{sec:non-embedded}

\begin{definition}\label{def:non-embedded} 
{\rm A quantum homogeneous $G$--space $(\calf_q(V),\varphi)$ 
is {\it embedded} if there is a right coideal subalgebra 
$\cala$ in $\funqg$ such that 
\begin{itemize}
\item[(i)] $\calf_q(V)$ is isomorphic to $\cala$ as an $\funqg$--comodule algebra, where the coaction on the latter is the restriction to $\cala$ of 
the comultiplication $\Delta$. 
\item[(ii)] The quotient $\mk$--module $\funqg/\cala$ is torsion free. 
\end{itemize}
}
\end{definition}

Condition (i) above is standard, see for example Section 11.6.1 in \cite{ks}. 
Condition (ii) ensures that the natural surjection of 
$\calf_q(V)$ onto its `classical limit' $\calf_q(V)/(q-1)\calf_q(V)$ can be  
identified with the restriction to $\cala$ of the natural 
surjection $\funqg\to\fung$.  

The following result shows that our $\funqorb$ is a non-embedded quantum homogeneous 
space for most $\xi\in P$. 
This indicates that the use of the interplay between $\lqm$ and $\funqm$ 
due to \cite{bm}, \cite{majid} combined with the `orbit map' of \cite{dfl} 
can not be replaced by an orbit map going directly from $\lqm$ to 
$\funqg$ (like in \cite{dm2}). 

\begin{proposition}\label{prop:embedded} 
Assume that for $\xi\in P$, one of the following holds: 
\begin{itemize}
\item[(i)] the size of the nilpotent Jordan block of $\xi$ is at least three; 
\item[(ii)] $n\geq 3$, and the size of the nilpotent Jordan block of $\xi$ is two; 
\item[(iii)] $\xi$ has at least three (different) non-zero eigenvalues. 
\end{itemize} 
Then $\funqorb$ is not an embedded quantum homogeneous space. 
\end{proposition}

\begin{proof} 
Assume in the contrary that for a $\xi$ satisfying one of the conditions 
(i), (ii), (iii), there exists a $\funqg$--comodule algebra injection 
$\iota: \funqorb\to\funqg$ such that $\funqg/\image(\iota)$ is torsion free 
as a $\mk$--module. Then we have a commutative diagram 
of algebra homomorphisms 
\[\begin{array}{ccccc} 
\funqorb &\stackrel{\iota}\longrightarrow &\funqg
&\stackrel{\varepsilon}\longrightarrow &\mk\\
\downarrow & &\downarrow & &\downarrow\\
\calf(\orbclxi)&\stackrel{\bar\iota}\longrightarrow &\calf(G)
&\stackrel{\bar\varepsilon}\longrightarrow&\mc
\end{array}
\] 
where $\varepsilon$, $\bar\varepsilon$ denote the counit map of the corresponding Hopf algebra, 
and the vertical arrows are the natural surjections with kernel generated by 
$q-1$. 
Since $\funqg/\image(\iota)$ is torsion free as a $\mk$--module, 
the map $\bar\iota$ is also injective. 
Recall that $\funqorb$ is defined as a quotient of $\lqm$. 
Keep the notation $l^i_j$ for the images of the generators of $\lqm$ in 
$\funqorb$, as well as in $\calf(\orbclxi)$. The fact that $\bar\iota$ is injective implies that the complex $n\times n$ matrix 
$\overline{B}=[\bar\varepsilon\bar\iota(l^i_j)]_{n\times n}$ belongs to the same $G$--orbit as $\xi$. On the other hand, this matrix is obtained by reducing 
modulo $(q-1)$ the $n\times n$ matrix 
$B=[\varepsilon\iota(l^i_j)]_{n\times n}$ 
with entries in $\mk$. Now $B$ satisfies the reflection equation 
(\ref{eq:reflection}). 
From the classification of solutions of the reflection equation 
given in \cite{mudrov} we see that it is impossible for a solution $B$ that 
$\overline{B}$ has the same Jordan normal form as $\xi$ 
(the paper \cite{mudrov} classifies complex solutions, but the proofs and the result can easily be extended to the case when $\mc$ is replaced by the domain 
$\mk$). This contradiction finishes the proof. 
\end{proof}


\section{Multiplicities of irreducible corepresentations}
\label{sec:multiplicities}

Let us recall that the corepresentation theory of the $\mc(q)$--Hopf algebra 
$\funkqg$ is completely analogous to its classical counterpart for $\fung$ 
(i.e. the representation theory of $G$ expressed in a dual language). 
Namely, $\funkqg$ is cosemisimple, and its irreducible corepresentations are indexed by $\Omega$, the set of dominant integral weights for $G$; 
see for example the book \cite{ks}. 
For $\lambda\in\Omega$ denote 
$\varphi_{\lambda}:V_{\lambda}\to V_{\lambda}\otimes_{\mc(q)}\funkqg$ 
the irreducible corepresentation indexed by $\lambda$. 
Here $V_{\lambda}$ has a decomposition into a direct sum of weight subspaces having the same dimension as in the classical case of $\fung$. 
Given $\lambda\in\Omega$ denote by $c(\lambda)$ the dimension of the 
trivial weight subspace in the dual of the irreducible 
$SL(n,\mc)$--representation with highest weight $\lambda$. 
It is proved in \cite{kostant} that for any 
$\xi\in M(n,\mc)$ whose orbit is maximal, the multiplicity of the 
irreducible $SL(n,\mc)$--module associated to $\lambda$ in the coordinate 
ring $\calf(\orbclxi)$ is $c(\lambda)$ for all $\lambda\in\Omega$. 
The fact that $\funqorb$ is an appropriate deformation of $\calf(\orbclxi)$ can be expressed in terms of corepresentation theory as follows. 
Extend scalars and set $\funkqorb=\mc(q)\otimes_{\mk}\funqorb$. 
Since $\funqorb$ is a flat $\mk$--module, the $\mc(q)$--algebra 
$\funkqorb$ has the same presentation in terms of generators and relations over $\mc(q)$ as the algebra $\funqorb$ over $\mk$. 
Now $\funkqorb$ is a right comodule algebra over $\funkqg$ and we have the 
following decomposition of this corepresentation. 

\begin{theorem}\label{thm:multiplicity} 
The corepresentation 
$\beta:\funkqorb\to\funkqorb\otimes_{\mc(q)}\funkqg$ 
for $\xi\in P$ 
decomposes into a direct sum of irreducible corepresentations as 
\[\beta\cong\sum^\oplus_{\lambda\in\Omega}c(\lambda)\varphi_{\lambda}.\]
\end{theorem}

\begin{proof} 
For $d=0,1,2,\ldots$, denote by $\funkqorb^d$ the $\mc(q)$--subspace spanned 
by the products of length $\leq d$ in the generators $l^i_j$ of 
$\funkqorb$. This is a finite dimensional subcomodule, containing the free 
$\mk$--submodule $\funqorb^d$. A weight subspace of $\funkqorb^d$ 
of dimension $r$ intersects $\funqorb^d$ in a free $\mk$--submodule of rank $r$. 
Under the natural surjection 
$\funqorb^d/(q-1)\funqorb^d\cong \calf(\orbclxi)^d$ 
this $\mk$--submodule is mapped onto an $r$--dimensional classical 
weight subspace belonging to the same weight of $G$. 
Therefore the dimensions of the weight subspaces in 
$\funkqorb^d$ agree with the dimensions of the corresponding weight subspaces 
in $\calf(\orbclxi)^d$. The multiplicities of the irreducible summands 
can be calculated from the weight multiplicities by the same rule for $\funkqg$ as for $\fung$, 
hence the statement follows from the corresponding classical result for 
$G$. 
\end{proof} 


\section{Quantum spheres}\label{sec:spheres}

Podle\'s \cite{p} introduced a family 
of quantum $2$--spheres, which became frequently 
cited examples in the literature on quantum homogeneous spaces. 
Here we point out that they can be obtained as coadjoint orbits in the 
case $n=2$. 
In this Section we change the setup and take $k=\mc$ as our base ring, 
so $q\in\mc^{\times}$ is a non-zero complex number. 
The algebra 
$\lcqm$ in the special case $n=2$ is generated by 
$l_{11}$, $l_{12}$, $l_{21}$, $l_{22}$, subject to the relations

\begin{eqnarray*}\label{eq:2rea} 
l_{22}l_{12}=q^2l_{12}l_{22};&\qquad
l_{11}l_{12}=l_{12}l_{11}+(q^{-2}-1)l_{12}l_{22};&\\ 
l_{11}l_{22}=l_{22}l_{11};&\qquad
l_{21}l_{12}=l_{12}l_{21}+(q^{-2}-1)l_{22}(l_{22}-l_{11});&\\ 
l_{21}l_{22}=q^2l_{22}l_{21};&\qquad
l_{21}l_{11}=l_{11}l_{21}+(q^{-2}-1)l_{22}l_{21}.&
\end{eqnarray*}

The basic coinvariants in $\lcqm$ are 
\[\Phi(\tau_1)=ql_{11}+q^{-1}l_{22}=\qtr(L)\] 
and  (after an easy calculation) 
\[\Phi(\tau_2)=l_{11}l_{22}-q^2l_{12}l_{21}
=(q+q^{-1})^{-1}(q\qtr(L)^2-q^2\qtr(L^2)).\] 

Classically a maximal $2\times 2$ coadjoint orbit is determined by specifying the 
values of the trace and the determinant. 
Motivated by this (or by the previous Sections), for a pair of complex parameters 
$(t,d)$ we consider the quotient algebra 
\[\lcqtd=\lcqm/(\qtr(L)=t,\Phi(\tau_2)=d).\] 

With the notation of \cite{ks}, for a pair of complex parameters 
$\alpha,\beta$, the Podle\'s sphere 
$\qsphere$ is generated by $x_{-1},x_0,x_1$, 
subject to the relations (71)--(74) on page 124 of \cite{ks}. 

\begin{proposition}\label{prop:podles} 
The $\funqg$--comodule algebra $\lcqtd$ is isomorphic to 
$\qsphere$, the Podle\'s sphere with parameters 
$\alpha=q^{-1}(q+q^{-1})^{-1/2}t$, 
$\beta=\alpha^2-(q^{-1}+q^{-3})d$. 
\end{proposition}

\begin{proof}  
Both algebras are defined in terms of generators and relations. 
Write $l_{ij}$ for the image in $\lcqtd$ of the corresponding generator of 
$\lcqm$. By the relation $\qtr(L)=t$ in $\lcqtd$ we can eliminate one of the generators. Thus $\lcqtd$ is generated by 
$il_{12},il_{21},(q+q^{-1})^{-1/2}(l_{11}-l_{22})$, where $i\in\mc$ is the imaginary unit with $i^2=-1$. 
Expressing the defining relations of $\lcqtd$ in terms of these three generators 
one gets the defining relations of $\qsphere$, so 
the map $il_{21}\mapsto x_{-1}$, $il_{12}\mapsto x_1$, 
$(q+q^{-1})^{-1/2}(l_{11}-l_{22})\mapsto x_0$ extends to an algebra isomorphism 
$\lcqtd\to\qsphere$. 
\end{proof}

\bigskip
\noindent{\bf Acknowledgement} 

The author thanks Pham Ngoc \'Anh, Rita Fioresi, and Tom Lenagan 
for conversations on the topic of this paper. 



\begin{thebibliography}{MMMM}

\bibitem{bfgp} P. Bonneau, M. Flato, M. Gerstenhaber, and 
G. Pinczon, 
The hidden group structure of quantum groups: strong duality, rigidity and 
preferred deformations, 
Comm. Math. Phys. 161 (1994), 125-156. 

\bibitem{bm} 
T. Brzezi\'nski and S. Majid, 
A class of bicovariant differential calculi on Hopf algebras, 
Lett. Math. Phys. 26 (1992), 67-78. 


\bibitem{em}
H. Cartan and S. Eilenberg, 
Homological Algebra, 
Princeton University Press, Princeton, 1956 (1999). 

\bibitem{dfl} 
M. Domokos, R. Fioresi and T. H. Lenagan, 
Orbits for the adjoint coaction on quantum matrices, 
J. Geometry and Physics 47 (2003), 447-468.

\bibitem{dl} 
M. Domokos and T. H. Lenagan, 
Conjugation coinvariants of quantum matrices, 
Bull. London Math. Soc. 35 (2003) 117-127. 

\bibitem{dl2}
M. Domokos and T. H. Lenagan, 
Weakly multiplicative coactions of quantized function algebras, 
J. Pure Appl. Alg. 183 (2003) 45-60. 

\bibitem{dm1} 
J. Donin and A. Mudrov, 
$U_q({\rm sl}(n))$--covariant quantization of symmetric coadjoint orbits via reflection equation algebra, 
Quantization, Poisson brackets and beyond (Manchester, 2001), 61--79, 
Contemp. Math., 315, 
Amer. Math. Soc., Providence, RI, 2002. 


\bibitem{dm2} 
J. Donin and A. Mudrov, 
Method of quantum characters in equivariant quantization, 
Comm. Math. Phys. 234 (2003), 533-555. 

\bibitem{dm3}
J. Donin and A. Mudrov, 
Explicit equivariant quantization on coadjoint orbits of 
${\rm GL}(n,\mathbb C)$, 
Lett. Math. Phys. 62 (2002), no. 1, 17-32. 

\bibitem{drinfeld}
V. G. Drinfeld, 
Quantum groups. 
Proceedings of the International Congress of Mathematicians, Vol. 1, 2 (Berkeley, Calif., 1986), 798--820, 
Amer. Math. Soc., Providence, RI, 1987. 

\bibitem{gl} 
K. R. Goodearl and T. H. Lenagan, 
Quantum determinantal ideals, 
Duke Math. J. 103 (2000), 165-190. 

\bibitem{gur-sap} 
D. Gurevich and P. Saponov, 
Quantum line bundles via Cayley-Hamilton identity, 
J. Phys. A: Math. Gen. 34 (2001), 4553-4569. 

\bibitem{hayashi}
T. Hayashi, 
Quantum deformations of classical groups, 
Publ. RIMS Kyoto Univ. 28 (1992), 57-81.

\bibitem{iop} A. Isaev, O. Ogievetsky, and P. Pyatov, 
Generalized Cayley-Hamilton-Newton identities, 
Czechoslovak J. Phys. 48 (1998), no. 11, 1369-1374.

\bibitem{jimbo}
Michio Jimbo, 
A $q$-analogue of $U({\mathfrak{gl}}(N+1))$, Hecke algebra, and the Yang-Baxter equation. 
Lett. Math. Phys. 11 (1986), 247-252.


\bibitem{kirillov}
A. A. Kirillov, 
Lectures on the Orbit Method, 
Amer. Math. Soc., Povidence RI, 2004. 

\bibitem{ks}
A. Klimyk and K. Schm\"udgen, 
Quantum Groups and Their Representations, 
Springer-Verlag, Berlin, Heidelberg, New York, 1997. 

\bibitem{kostant}
B. Kostant, 
Lie group representations on polynomial rings, 
Amer. J. Math. 85 (1963), 327-404. 

\bibitem{kraft} 
H. Kraft, 
Geometrische Methoden in der Invariantentheorie, (German), 
Vieweg Verlag, Braunschweig-Wiesbaden, 1985. 

\bibitem{kulish-sk} 
P. P. Kulish and E. K. Sklyanin, 
Algebraic structures related to reflection equations, 
J. Phys. A: Math. Gen. 25 (1992), 5963-5975. 

\bibitem{majid} 
S. Majid, 
Foundations of Quantum Group Theory, 
Cambridge Univ. Press, Cambridge, 1995. 

\bibitem{majid2}
S. Majid, 
Examples of braided groups and braided matrices, 
J. Math. Phys. 31 (1991), 3246-3253. 

\bibitem{mudrov}
A. Mudrov, 
Characters of $U_q(gl(n))$--reflection equation algebra, 
Lett. Math. Phys. 60 (2002), 283-291. 

\bibitem{nym} 
M. Noumi, H. Yamada, and K. Mimachi, Finite
dimensional representations of the quantum group $GL_q(n;\mc)$ and
the zonal spherical functions on $U_q(n-1)\backslash U_q(n)$, Japanese
J. Math. {\bf 19} (1993),  31-80.

\bibitem{p} 
P. Podle\'s, 
Quantum spheres, 
Lett. Math. Phys. {\bf 14} (1987), 193-202. 

\bibitem{rtf}
N. Yu. Reshetikhin, L. A. Takhtadzhyan, and L. D. Faddeev,
Quantization of Lie groups and Lie algebras, (Russian), 
Algebra i Analiz 1 (1989), 178--206.


\end{thebibliography}
\end{document}